\documentclass[a4paper, 11pt]{article}
\usepackage{xcolor}
\usepackage[english]{babel}
\usepackage{amssymb,amsmath,graphicx}
\usepackage{amsthm}
\usepackage{mathtools}
\usepackage[ruled, vlined, linesnumbered]{algorithm2e}
\usepackage{setspace}			
\usepackage[round]{natbib}
\usepackage{pdflscape}			
\usepackage{ushort}
\usepackage{microtype}
\usepackage{tabularx}
\usepackage{enumitem}

\usepackage[margin=1.5in]{geometry}
\newcolumntype{C}{>{\centering\arraybackslash}p{2.5em}}

\begin{document}

\begin{center}

\vspace*{-1.5cm}

\begin{LARGE}
Heuristic Rectangle Splitting:
\linebreak
Leveraging single-objective heuristics \\ to efficiently solve multi-objective problems
\end{LARGE}

\vspace*{1.00cm}


\textbf{P.~Matl*} \\
University of Vienna, Austria \\
piotr.matl@univie.ac.at \\

\vspace*{0.2cm}
\textbf{R.F.~Hartl} \\
University of Vienna, Austria \\
richard.hartl@univie.ac.at \\

\vspace*{0.2cm}
\textbf{T.~Vidal} \\
Pontifical Catholic University of Rio de Janeiro, Brazil  \\
vidalt@inf.puc-rio.br  \\


\vspace*{0.7cm}




\thispagestyle{empty}

\begin{quotation}

Real-life problems are often characterized by conflicting optimization objectives. Consequently, there has been a growing interest not only in multi-objective models, but also in specialized multi-objective metaheuristics for solving those models. A wide variety of methods, e.g.~NSGA-II, SPEA, IBEA, scatter search, Pareto local search, and many others, have thus been proposed over the years. Yet in principle, multi-objective problems can be efficiently solved with existing tailored single-objective solvers -- this is the central idea behind the well-known $\epsilon$-constraint method (ECM). Despite its theoretical properties and conceptual simplicity, the $\epsilon$-constraint method has been largely ignored in the domain of heuristics and remains associated mostly with exact algorithms. In this article we dispel these preconceptions and demonstrate that the $\epsilon$-constraint framework can be a highly effective way to directly leverage the existing research on single-objective optimization for solving multi-objective problems.

We propose an improved version of the classical ECM adapted to the challenges and requirements specific to heuristic search. The resulting framework is implemented with an existing state-of-the-art single-objective solver for the Capacitated Vehicle Routing Problem (CVRP) and tested on the VRP with Route Balancing (VRPRB). Based on an extensive computational study, we show the added value of our adaptations compared to the classical ECM, and demonstrate that our simple $\epsilon$-constraint algorithm significantly outperforms the current state-of-the-art multi-objective metaheuristics with respect to multiple quality metrics. We conclude with a discussion of relevant success factors and promising directions for further research.


\bigskip
\newline * Corresponding author

\end{quotation}
\end{center}


\newpage
\onehalfspacing
\pagenumbering{arabic}

\section{Introduction}
\label{1-intro}

In the decades following the introduction of the vehicle routing problem (VRP) by Dantzig and Ramser in 1959, much of the research on the VRP and its variants has focused on models and methods in which the objective is to minimize cost or distance \citep{toth2014}, leading to significant advances on single-objective (SO) heuristic algorithms. In recent years, an increasing number of VRP variants have been presented that also account for other relevant objectives which conflict with direct cost minimization. These include service consistency \citep{kovacs2014}, workload balancing \citep{matl2016}, safety and security \citep{huang2004, samanlioglu2013}, service quality \citep{park2010}, as well as emissions and environmental externalities \citep{lin2014}. As a result, there has been increasing interest also in multi-objective (MO) VRP algorithms, which aim to determine a set of compromise solutions. Decision makers can then gain greater insight into the trade-offs made between conflicting objectives.

Intuitively, one may consider it self-evident that MO problems require dedicated solution approaches. After all, SO algorithms aim to find \textit{one} optimal solution, whereas in MO problems, no single optimal solution exists in general. And indeed, just as research on SO optimization has led to the development of a variety of metaheuristics, the research on MO optimization has resulted in an even larger set of sophisticated methods, such as NSGA-II, IBEA, variants of scatter search and evolutionary algorithms, Pareto local search, path re-linking, among many~others.

Yet in principle, the complete Pareto set for discrete MO problems can be found using only one solver optimizing one of the respective objectives: this is the fundamental property underlying the classical $\epsilon$-constraint method. By relegating all other objectives into the constraint space or into penalty terms in the primary objective function, it is possible to generate the entire set of Pareto-optimal solutions by iteratively solving modified SO sub-problems. 

Although the $\epsilon$-constraint method is well-established within the domain of exact optimization, it has received very limited attention within the field heuristics. In a survey of MO VRP models and methods, \cite{jozefowiez2008} conclude that the vast majority of proposed algorithms use weighted aggregations of objectives or, more commonly, some form of genetic algorithm. This is surprising, since in principle only minor adjustments are needed to embed an existing tailored SO solver into an $\epsilon$-constraint framework. And this in turn allows to directly leverage the decades of research on SO heuristics.

\newpage
In this article, we revisit the classical $\epsilon$-constraint framework from a heuristic perspective. We propose a variation of the classical approach that is better suited to the characteristics and demands of heuristic search, and detail how an existing tailored state-of-the-art CVRP solver can be integrated into our algorithm with only limited changes. An in-depth computational analysis reveals the added value of our modifications compared to the classical $\epsilon$-constraint method. Finally, we demonstrate through an extensive empirical study that this simple and generic approach significantly outperforms the current state-of-the-art methods for the VRPRB in terms of solution quality, reliability, and computational effort.






\section{Literature Review}
\label{2-litreview}

A general overview of MO VRP algorithms published up to 2008 can be found in the survey by \cite{jozefowiez2008}. In this review we focus our attention specifically on the heuristics presented thus far for the VRP with Route Balancing (VRPRB). The VRPRB extends the classical CVRP with a second objective aimed at reducing the imbalance (inequity) between the vehicle workloads. This is measured as the range of their tour lengths, or more formally as $Z=\max_i d_i - \min_i d_i$, where $d_i$ represents the distance or duration of tour $i$. A mathematical formulation can be found in \cite{oyola2014}, but some conventions vary from paper to paper (e.g.~the fleet size and possible 2-optimality of tours).

The VRPRB was first presented in \cite{jozefowiez2002}. A parallel genetic algorithm with tabu search was proposed as a solution procedure. This line of research was later extended in \cite{jozefowiez2006} with a parallel version of NSGA-II combined with enhanced strategies for managing population diversity. This, in turn, was followed by the introduction of target-aiming Pareto search, a hybridization of NSGA-II, tabu search, and goal programming \citep{jozefowiez2007}.

Following these initial contributions, alternative solution procedures have been proposed by other authors. \cite{pasia2007a} describe a population-based local search method based on the concept of Pareto local search. The authors show that generating initial solutions with a randomized Clarke-Wright savings algorithm helps to improve significantly the performance of the method, and a comparative study demonstrates that their algorithm outperforms the one proposed by \cite{jozefowiez2006}. These results are further improved in \cite{pasia2007b} with extensions based on Pareto ant colony optimization, which adds a learning layer to the previous method. Adaptive memory mechanisms are explored also by \cite{borgulya2008}, who replaces the standard recombination step of an evolutionary algorithm with an adaptive mutation operator. Finally, \cite{jozefowiez2009} revisit the VRPRB and present a classical evolutionary algorithm extended with parallel search and elitist diversification management. Computational experiments demonstrate the positive impact of both new mechanisms on the quality of the generated non-dominated solution sets.

Unlike most earlier approaches, more recent contributions on the VRPRB do not propose MO evolutionary algorithms. \cite{oyola2014} present a bi-objective GRASP heuristic (GRASP-ASP). Rather than beginning each iteration from an empty solution, a starting solution (``advanced starting point'') is formed from the common elements of two previously found non-dominated solutions. This partial solution is then extended with a bi-objective GRASP procedure which considers for each customer insertion the impact on both objective functions and the resulting solution's Pareto rank. A bi-objective local search procedure is applied to each complete solution before evaluating its inclusion in the archive of non-dominated solutions. On small instances with up to 100 customers, GRASP-ASP matches or outperforms the best bounds found by an exact method based on a weighted sum objective. The authors also reimplemented a sequential version of the algorithm proposed in \cite{jozefowiez2009}, yielding slightly better solution quality, but at the expense of considerably longer computation times.

The most recent article on the VRPRB \citep{lacomme2015} introduces a multi-start decoder-based method (MSSPR) which alternates between direct and indirect solution spaces. Non-dominated sets of solutions are extracted from TSP giant tours through a bi-objective \textsc{Split} procedure \citep{prins2004, vidal2016}. Similar to GRASP-ASP, a bi-objective local search is applied to all solutions found. The resulting VRPRB solutions serve as starting points for further search based on path re-linking, and a multi-start strategy is applied for diversification. Although an in-depth comparative study with previous methods was not possible due to a lack of detailed computational results for previous algorithms, the MSSPR of \cite{lacomme2015} is shown to be reasonably competitive on the cost objective when compared to the best known SO solutions.

\section{The $\epsilon$-Constraint Framework with Heuristics}
\label{3-0}

In what follows, we consider an optimization problem with two minimization objectives: $f_1$ and $f_2$ (maximization objectives can be appropriately transformed). In general there is no single optimal solution minimizing both objectives simultaneously, and all solutions imply some degree of compromise. A trade-off solution $\mathbf{x}$ in the set $\mathbf{X}$ of feasible solutions is Pareto-\textit{optimal} if it is impossible to improve either objective without worsening the other. Formally: there exists no other solution $\mathbf{x^\ast} \in \mathbf{X}$ such that $f_1(\mathbf{x^\ast}) \leq f_1(\mathbf{x})$ and $f_2(\mathbf{x^\ast}) \leq f_2(\mathbf{x})$, with at least one strict inequality. In order to more clearly distinguish heuristic solutions that approximate the optimal Pareto set, we will say that a solution is merely \textit{non-dominated} if the above conditions hold only for a given subset of $\mathbf{X}$. Our aim is to identify the set of Pareto-optimal trade-off solutions, or a suitable approximation of mutually non-dominated solutions.

In principle, MO optimization problems can be efficiently solved by embedding a single-objective solver within the well-known $\epsilon$-constraint method (ECM). By iteratively constraining at least one objective to strictly improve compared to all previously found solutions, all Pareto-optimal solutions are guaranteed to be identified. A formal description of the ECM as well as an efficient generalization to more than two objectives is given in \cite{laumanns2006}.

Although the ECM is commonly used in the field of exact optimization, its application has remained largely unexplored in the domain of heuristics. Yet the framework offers several appealing properties:

\paragraph{Advantages of the Classical $\epsilon$-Constraint Method}
\begin{itemize}
	\item By design, the ECM ensures finding in every iteration a solution which is non-dominated with respect to all previously identified solutions. As a result, each iteration strictly improves the quality of the current Pareto set, regardless of whether the solver is exact or heuristic.
	
	\item The ECM is well-suited to exploit the possible similarities between consecutive non-dominated solutions. With exact methods, \cite{boland2015} demonstrate the value of re-using high-quality solutions as starting points for subsequent sub-problems. This translates in a similar way to heuristics, in which re-using previous local optima, adaptive parameters, etc.~can greatly speed up and improve the search for other non-dominated solutions.
	
	\item The performance of the ECM, in terms of both solution quality and computational effort, is tied directly to the quality and efficiency of its underlying SO solver. Hence all problem-specific aspects are contained within a single sub-procedure and free from dependencies with higher level search strategies -- this cannot be said for many MO heuristics.
	
	\item Finally, the ECM is generic, easy to understand, and straightforward to implement. It also requires the specification of only a single parameter -- the eponymous $\epsilon$ step size.
\end{itemize}

Despite these advantages, some difficulties can arise when using the $\epsilon$-constraint method with a heuristic solver. These issues may explain why it has not been as readily adopted with heuristics as with exact approaches.

\paragraph{Challenges with Heuristics}
\begin{itemize}

	\item Without a strict optimality guarantee, it is possible that previously found heuristic solutions become dominated, even late into the search. Although this still constitutes a strict improvement of the Pareto front, this can lead to gaps in the generated solution set, such that some regions of the objective space are not approximated by any solution. As a result, the classical ECM might not be appropriate, since it never returns to previously examined areas of the objective space.
	
	\item It is desirable for algorithms to converge to good approximations quickly, especially in a heuristic context. Yet the classical ECM has very poor ``anytime'' behavior, leaving entire regions of the objective space without even a single solution if terminated prematurely.
	
	
	\item Finally, setting the $\epsilon$ parameter is far from trivial, and the ideal setting is instance-dependent. With the minimal value, the computational effort can become exponential with the number of Pareto-optimal solutions. On the other hand, larger values, whether in absolute or percentage terms, can skip significant parts of the Pareto set. The extent of this behavior depends not only on the quality of lower bounds for the constrained objective, but also on the distribution of solutions in the objective space -- yet this information is unknown ex ante.
	
\end{itemize}

Our assertion is that the essential properties and advantages of $\epsilon$-constraint frameworks can be exploited also when the underlying solver is a heuristic, and we demonstrate that the potential difficulties can be mitigated with appropriate adaptations to the framework.

In Section \ref{3-1} we recall an existing  adaptation of the classical ECM which resolves some of the challenges above, and then propose improvements tailored for heuristic solvers. To demonstrate its performance, we test our approach on the VRPB using a state-of-the-art single-objective VRP solver:~the Hybrid Genetic Search (HGS) of \cite{vidal2012}. We recall in Section \ref{3-2} the general behaviour of the underlying HGS, and explain the adjustments required to embed it within an $\epsilon$-constraint-based framework in Section \ref{3-3}. The results of our computational experiments are reported and discussed in the sections thereafter.


\subsection{Heuristic Rectangle Splitting}
\label{3-1}

Since the set of non-dominated solutions to a MO problem can be of exponential size, we focus our attention on generating a representative subset within bounded computational time. Because non-dominated solutions can be distributed arbitrarily in the objective space, an effective algorithm must be capable of adapting to the true distribution of the trade-off solutions and directing computational effort accordingly.

\paragraph{Underlying Principles}
Algorithms for generating Pareto set representations are often based around splitting the objective space according to some rule. The Box Algorithm proposed by \cite{hamacher2007} is a prototypical template. The method begins by solving two lexicographic optimization problems to identify the two extreme solutions with minimal $f_1$ and minimal $f_2$, respectively. These two solutions represent a rectangle containing the entire Pareto set. By appropriately setting the $\epsilon$-constraint, the objective space represented by the rectangle can be temporarily split into feasible and infeasible halves. The solution to the corresponding sub-problem allows to discard areas in the objective space which cannot contain non-dominated solutions, hence splitting the original rectangle into at most two smaller ones. This procedure is then repeated by selecting in each iteration the largest remaining rectangle, in effect performing a type of binary search in the objective space.

In principle, the Box Algorithm resolves some of the difficulties and disadvantages inherent to the classical ECM. By always splitting the largest remaining box, the method quickly converges to a representative approximation. It allocates computational effort to those areas in objective space which are least represented by the current solution set and which are most likely to contain further non-dominated solutions. Finally, it eliminates the $\epsilon$ parameter entirely, avoiding the need to select a ``good'' $\epsilon$ ex ante.

However, the original article by \cite{hamacher2007} did not report empirical tests, and the Box Algorithm can be further improved and generalized. In the context of exact approaches, \cite{boland2015} demonstrate that solving two appropriately chosen lexicographic sub-problems per rectangle increases the rate of convergence and that warm-starting subsequent sub-problems with previously found solutions (``\textit{solution harvesting}'') reduces computational effort even further. In the context of heuristics, it is furthermore necessary to account for the lack of a strict optimality guarantee when updating and discarding areas in the objective space during the splitting procedure. In addition, the computation of both extreme points generally requires \textit{two} dedicated solution methods. State-of-the-art single-objective algorithms, whether exact or heuristic, exploit properties of the objective function or the problem structure. Simply changing the objective function is unlikely to be efficient with such methods, especially when faced with atypical, complex, and/or non-linear objectives. 

In the following, we propose a \textbf{Heuristic Rectangle Splitting (HRS)} algorithm which resolves these issues without compromising the principles underpinning the exact approaches. We rely on only one state-of-the-art solver for one of the objectives.

\begin{landscape}
\begin{figure} [!t]
	\centering
		\includegraphics[width=1.65\textwidth]{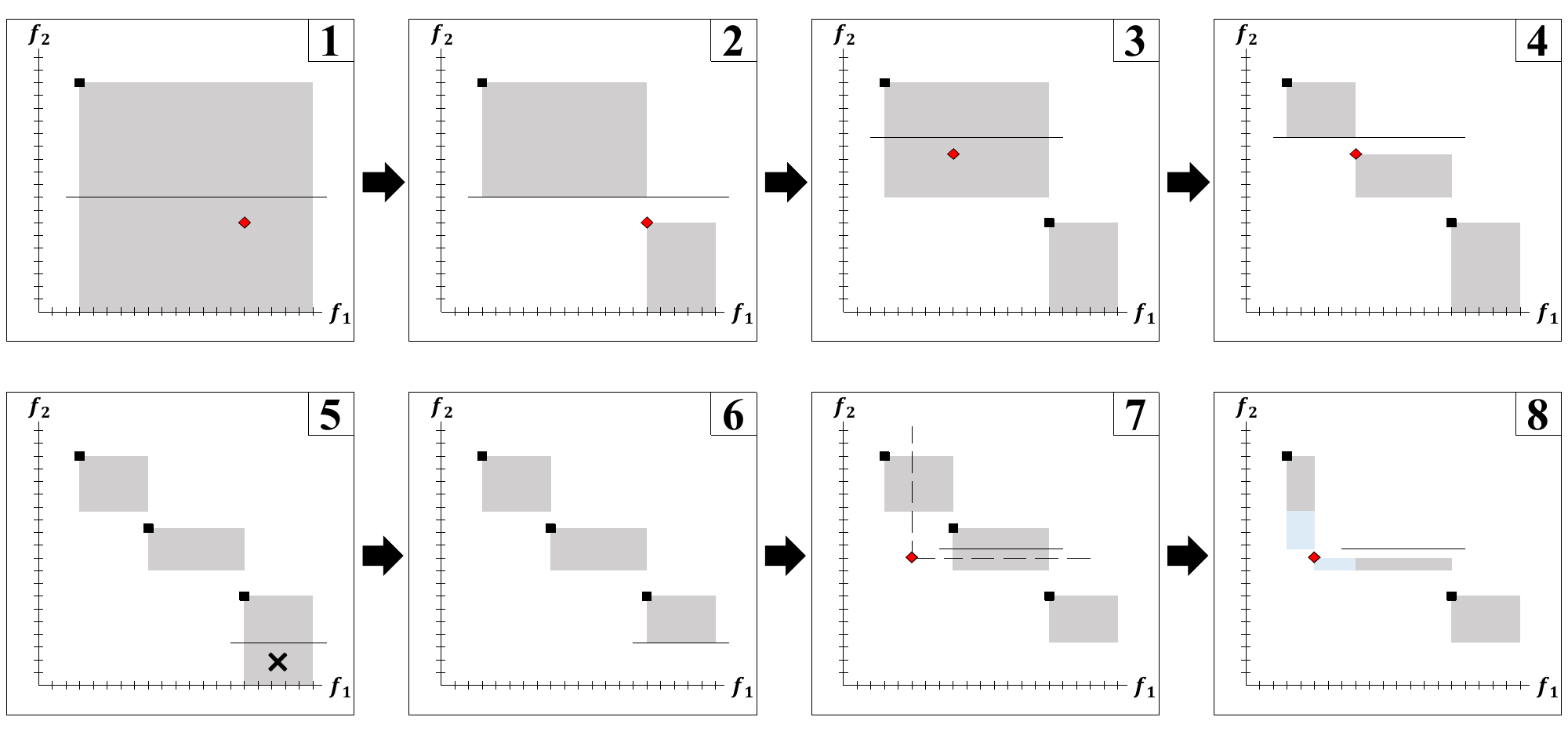}
	\caption{Visualization of rectangle splitting and updating. Black squares represent solutions, shaded areas delimit their associated rectangles. The $\epsilon$-constraints are indicated with a horizontal line, and the solutions to the corresponding sub-problems are marked with a diamond (shaded red in the electronic version of this article). The lightly shaded areas in step 8 indicate segments which would have been added according to the updating procedure of the original Box Algorithm.}
	\label{HRSviz}
\end{figure}
\end{landscape}

\paragraph{Notations}
For convenience, we will denote with $R(z^1,z^2)$ the rectangle formed by the points $(z^1_1, z^1_2)$ and $(z^2_1,z^2_2)$ in the objective space, such that $z^1_1 \leq z^2_1$ and $z^1_2 \geq z^2_2$ (i.e.~point $z^1$ is the ``upper-left'' vertex and point $z^2$ the ``lower-right'').
The set of non-dominated solutions is maintained in an archive $A$, the set of rectangles in an archive $\mathcal{R}$. With \textsc{optimize}($f_1, c$) we denote an algorithm which minimizes objective $f_1$ subject to the constraint $f_2 \leq c$. With $f^{min}_2$ we refer to any valid lower bound on objective $f_2$, and with $f^{max}_1$ any upper bound on objective $f_1$ (note that both may be trivial, e.g.~0 and a big $M$, respectively). Algorithm \ref{HRS} describes the proposed method in numerical detail, and Figure \ref{HRSviz} illustrates the procedure. In the following, we describe each main element of the approach: the management of extreme points, and the rules to split and update the boxes.


\paragraph{Extreme Points}
Rather than relying on the availability or implementation of a second state-of-the-art solver for objective $f_2$, we will iteratively use the \textsc{optimize} solver to progressively improve an estimate of the $f_2$ extreme point while identifying other non-dominated solutions. Owing to the binary character of the splitting procedure, only a logarithmic number of sub-problems will be solved to reach the $f_2$ extreme point.
Depending on the quality of the lower bound $f^{min}_2$, as many as all of these sub-problems can yield non-dominated solutions that contribute to the Pareto set representation.

The algorithm begins by using the \textsc{optimize} procedure to minimize objective $f_1$ without any constraint on $f_2$. This yields the ``upper-left'' extreme point. The two bounds $f^{max}_1$ and $f^{min}_2$ correspond to the ``lower-right'' extreme point. Together, these two points delimit the initial rectangle. The approximation of the true ``lower-right'' extreme point will be progressively improved as the algorithm proceeds to split unexamined regions of the objective space and discard those which become redundant.


\begin{algorithm} [!t]
\DontPrintSemicolon
\KwIn{termination criterion, $f^{max}_1$, $f^{min}_2$}
\KwOut{Archive of non-dominated solutions $A$}
\SetKwBlock{Begin}{Initialize}{end Initialize}
$x \leftarrow$ \textsc{optimize}($f_1, \infty$)\;
	$A \leftarrow \left\lbrace x \right\rbrace, \ z^1 \leftarrow \left( x_1, x_2 \right), \ z^2 \leftarrow \left( f^{max}_1,f^{min}_2 \right), \ \mathcal{R} \leftarrow \left\lbrace R \left( z^1, z^2 \right) \right\rbrace$\;
	\Repeat {termination criterion reached}
	{
		Find $R(y^1,y^2) \in \mathcal{R}$ with maximal area\;
		$c \leftarrow \frac{1}{2}(y^1_2 + y^2_2)$\;
		$x \leftarrow$ \textsc{optimize}($f_1, c$)\;
		
		\If{$x$ \emph{is infeasible} $\mathbf{or}$ $x$ \emph{is dominated}}
		{
			$y^2 \leftarrow (y^2_1, c)$
		}
		\Else
		{
			Update archive $A$ with $x$\;
			
			\If{\emph{there exists} $\{R'(z^1,z^2) \in \mathcal{R} \mid z^1_1 < x_1 \leq z^2_1\}$}
			{
				$z^\ast \leftarrow (x_1,\max\{z^2_2, c\})$\;
				Add $R(z^1, z^\ast)$ to $\mathcal{R}$\;
			}
			\If{\emph{there exists} $\{R''(z^1,z^2) \in \mathcal{R} \mid z^1_2 \geq x_2 \geq z^2_2\}$}
			{
				$z^\ast \leftarrow (\max\{x_1, z^1_1\}, x_2)$\;
				
				Add $R(z^\ast, z^2)$ to $\mathcal{R}$\;
			}
			Delete $R'$ and $R''$ if they exist\;
			\ForAll{$R(z^1,z^2) \in \mathcal{R}$ \emph{with} $z^1$ \emph{dominated by} $x$}
			{
				Remove $R$ from $\mathcal{R}$ \;
			}
		}
		
			
			
	}
	
	\Return{$A$}
\caption{Heuristic Rectangle Splitting}
\label{HRS}
\end{algorithm}

\paragraph{Box Splitting and Update}
As in the original Box Algorithm, each iteration proceeds by setting the $\epsilon$-constraint such that the largest remaining rectangle is split in half and then solving the corresponding sub-problem (Figure \ref{HRSviz}-1). From the definition of Pareto-efficiency, it follows directly that the region dominated by the resulting solution $x$ can be removed from further consideration, regardless of whether the \textsc{optimize} solver is exact or heuristic. If the solver is exact, then the objective space defined by $f_1 \leq x_1$ and $f_2 \leq c$ can likewise be discarded with absolute certainty, as any solutions in this region would contradict the optimality of the solver \citep{hamacher2007}. In the absence of this strict optimality guarantee when \textsc{optimize} is a heuristic, discarding this second region follows a probabilistic interpretation of the previous argument, based directly on the quality and reliability of the heuristic. The resulting rectangles are depicted in Figure \ref{HRSviz}-2. 

The algorithm continues by iteratively selecting and splitting in half the largest remaining rectangle (Figure \ref{HRSviz}-3). The corresponding sub-problem typically yields a solution that lies within the boundaries of the selected rectangle, resulting in two smaller rectangles as shown in Figure \ref{HRSviz}-4.

Exceptions to this behaviour can arise if the initial ``lower-right'' extreme point consists of objective function bounds instead of the Pareto-optimal solution minimizing $f_2$. Under those conditions, it may occur that a sub-problem is infeasible and does not yield a non-dominated solution (Figure \ref{HRSviz}-5). In such cases, the archive remains unchanged, the explored half of the rectangle is completely discarded, and no new rectangle is created (Figure \ref{HRSviz}-6). In this way, also the rectangle which estimates the ``lower-right'' extreme point is successively reduced in half, and so the search quickly converges to a close approximation of the true extreme solution minimizing $f_2$.

If the \textsc{optimize} solver is a heuristic, then the lack of an optimality guarantee can lead to further exceptions for the original Box Algorithm. First, it is possible that newly identified solutions can lie outside the boundaries of the rectangle which was split in a given iteration. As a result, more than one box may need to be updated. Second, it is possible that a solution dominates previously identified solutions, or in fact entire rectangles in more extreme cases. These situations can occur simultaneously, as depicted in Figure \ref{HRSviz}-7.

In these cases, the rectangle update procedure of the original Box Algorithm leads to the \textit{enlargement} of existing rectangles. This is depicted by the lighter shaded regions in Figure \ref{HRSviz}-8. As this can potentially cause a long process of cycling sub-problems between the two enlarged rectangles, we adapt the original updating rules to prevent enlargements.


\paragraph{Summary}
We propose to extend the exact $\epsilon$-constraint-based Box Algorithm of \cite{hamacher2007} for use with heuristics. Our Heuristic Rectangle Splitting approach accounts for the lack of strict optimality guarantees with heuristic solvers, and requires only one such solver optimizing only one of the objectives. Since previously found heuristic solutions can be dominated in later stages of the search, we also propose alternative updating procedures for managing the regions of the objective space which remain promising. By performing a type of binary search in the objective space, our algorithm rapidly converges to a representative approximation of the true Pareto set. By always searching in the largest unexplored region, it adapts to the true distribution of the solutions in the objective space and allocates computational effort to those areas which are least represented at any given time. Finally, the problematic $\epsilon$ parameter is eliminated and can be replaced by a typical termination criterion such as the number of sub-problems solved.

Overall, the proposed HRS framework is a generic and flexible way to use existing single-objective heuristics to solve multi-objective problems. It exploits the inherent advantages of the classical ECM in a heuristic context while mitigating its main drawbacks. In the following sections, we describe the implementation of HRS with a state-of-the-art VRP solver.
\subsection{Single-Objective Hybrid Genetic Search}
\label{3-2}

The Hybrid Genetic Search (HGS) of \cite{vidal2012} attains state-of-the-art results on the CVRP in a single-objective setting, making it suitable for use as the underlying solver in our HRS framework. In the following, we briefly recall the main ideas behind HGS. For a more detailed description of the method we refer the reader to the original article and its e-companion. 

\paragraph{Evolution} HGS evolves a population of solutions, divided into feasible and infeasible sub-populations. Penalty terms for violated constraints are included in the objective function in order to allow for the exploration of infeasible solutions, while still biasing the search toward feasibility. Pairs of solutions for crossover are selected through binary tournament selection, whereby each solution's contribution to population diversity is included in its fitness score in addition to its objective function value (biased fitness). New candidate solutions, encoded as giant tours, are generated by applying the standard order crossover (OX) to a selected pair of parent solutions, and the resulting giant tour is segmented optimally into a VRP solution by means of the commonly used \textsc{Split} procedure \citep{prins2004, vidal2016}.

\paragraph{Local Search} All candidate solutions undergo a local improvement procedure consisting of nine classical intra and inter-tour moves (\textsc{Move}, \textsc{Swap}, \textsc{2-opt}, \textsc{2-opt*}, and their generalizations). Moves are evaluated in random order and immediately accepted if an improvement is found. Since moves involving very distant customers are unlikely to lead to improvements, the local search is restricted to moves involving nearby customers within a granularity threshold \citep{toth2003}. Once no further improving moves are found, the feasibility of the candidate solution is examined: if it is not feasible, an attempt is made to repair the solution by repeating the local search with temporarily increased penalty parameters. Either way, the resulting solution is added to the corresponding feasible or infeasible sub-population.

\paragraph{Population Management} Since the population size increases in each iteration, it eventually becomes necessary to discard some solutions to promote elitism. Whenever a sub-population reaches a given maximum size, the solutions are ranked according to their biased fitness and those with the worst biased fitness are removed until the sub-population reaches a given minimum size. Finally, the penalty parameters are adjusted adaptively during the search in order to favor the generation of feasible solutions: the more solutions are infeasible with respect to a particular constraint, the higher the penalty term for this violation, and vice versa.


\bigskip
Although HGS evolves a set of solutions, we emphasize that it is a single-objective solver and the genetic elements of the algorithm are \textit{not} designed to approximate a non-dominated set. In the context of HGS, population diversity does not relate to any other \textit{objective} functions, but rather to differences in the \textit{solution} space in terms of the edges selected in each VRP solution (broken-pairs distance). As a result, HGS promotes the diversification of the search in a systematic way but without neglecting the singular optimization objective.

\subsection{Integrating HGS into an $\epsilon$-Constraint Framework}
\label{3-3}


In general, an $\epsilon$-constraint can be handled either with a hard constraint on the secondary objective, or by introducing a corresponding penalty term to the primary objective. Since HGS already uses penalty functions to explore infeasible solutions, we simply introduce an additional penalty expression: the range of tour lengths in excess of the $\epsilon$-constraint threshold, multiplied with a penalty coefficient for imbalance. This penalty coefficient is adjusted in the same manner as the penalties for maximum route duration and vehicle capacity. A list of sorted route lengths is maintained for each solution in order to compute in constant time the change in the imbalance penalty during move evalutations.


The \textsc{Split} procedure remains largely unchanged. We use the limited fleet version of this algorithm, documented in \cite{prins2009}, and limit the maximum number of unused vehicles to one. Such a limit is necessary because classical local search moves cannot introduce customer visits into more than one empty tour simultaneously. As such, from an incumbent solution with two empty routes, it would be impossible to raise the minimum distance beyond $0$ in a single move, hindering progress towards better solutions in terms of the range equity objective. Moreover, we have compared this classical \textsc{Split} algorithm with a more sophisticated version that would be exact in the presence of the balance penalties. Using the exact version was computationally intensive but without effect on solution quality. This observation is in line with those of \cite{boudia2007} and \cite{lacomme2015}. A sub-optimal but faster \textsc{Split} algorithm is usually not a hindrance when embedded in a global search framework, because the task of improving a given solution is fully assumed by the local search. The crossover and \textsc{Split} procedures mainly serve to exploit past search information to provide meaningful starting solutions for further improvement.

\bigskip
The remaining adaptations concern how the algorithm behaves after the $\epsilon$-constraint is adjusted. Rather than re-initializing HGS for each new sub-problem, we save and re-use the evolved populations and adapted penalties (i.e.~the solver search state) from one sub-problem to the next. We associate with each non-dominated solution in the archive the search state at which the solution was retrieved, and when splitting rectangles we warm-start from the search state of the nearest non-dominated solution. Since neighboring non-dominated solutions tend to exhibit similar structure, this strategy significantly reduces the convergence time of HGS. For this reason, a comparatively large number of iterations without improvement ($I_0$) is allotted to the initial sub-problem of identifying the cost-optimum (as in the context of pure single-objective optimization), and a shorter limit ($I_{n}$) is used for all subsequent sub-problems, in order to efficiently explore the Pareto front. 

Since HGS performs only a brief search during each sub-problem, populations are diversified by introducing new random solutions at the start of each sub-problem (rather than after a fixed $I_{div}$ iterations without improvement as in the original method). Solutions from previous populations which become infeasible after the $\epsilon$-constraint is tightened are subjected to the repair procedure. Finally, HGS returns for each $\epsilon$-constraint sub-problem all feasible solutions in its populations (\textit{solution harvesting}). Although these other solutions are evolved only to diversify the search for the sub-problem's cost-optimum, there is a chance that some might be non-dominated, and hence all are examined for inclusion in the archive. However, the rectangles for the HRS framework are updated only based on the respective cost-optimum.

We will refer to the resulting algorithm as \textbf{Iterated HGS (IHGS)} in order to distinguish it from the more generic HRS framework. The pseudocode remains as in Algorithm \ref{HRS}, with the \textsc{optimize} sub-procedure replaced by the HGS described above.

\section{Computational Experiments}
\label{5-results}



\subsection{Parameter Settings}
\label{4-1}

Although the HGS parameters have remained largely unchanged from those used by \cite{vidal2012}, for clarity we (re)state all the relevant parameter values in Table \ref{params}. The impact of these parameters was evaluated during preliminary analyses. We observed that most of these values, which had been calibrated by \cite{vidal2012} for the single-objective CVRP, remain effective also for optimizing the single-objective sub-problems within our $\epsilon$-constraint framework.

Only two HGS parameter values -- $g$ and $\xi^{REF}$ -- were adapted to the bi-objective setting. Since the structure of highly balanced solutions is likely to differ from cost-minimizing solutions, the granularity parameter $g$ of the HGS local search has been doubled in order to explore more diverse neighbor solutions. Similarly, $\xi^{REF}$ has been doubled to allow for more feasible solutions and thus account for the additional source of infeasibility due to the balance constraint.

For the HRS framework, we specified a simple termination criterion of $n_{max}$ sub-problems. This is an instance and problem-independent CPU budget which ensures that small Pareto sets are explored exhaustively, while large sets are reasonably approximated with limited computational effort.

\begin{table}[h]
\centering
\resizebox{\textwidth}{!}{%
\begin{tabular}{rlc}
\hline
\multicolumn{2}{l}{\textbf{Parameter}}                                                  			& \textbf{Value} \\ \hline
$\mu$							& minimum population size (maximum is $\mu+\lambda$)                        & 25			\\
$\lambda$					& number of offspring per generation                                     		& 40			\\
$el$							& number of elite solutions kept in the population                       		& 10			\\
$\xi^{REF}$				& target proportion of feasible solutions in the population              		& 0.4			\\
$n_{close}$				& number of closest solutions considered for diversity evaluation           & 5 			\\
$p_{rep}$					& repair rate                                                            		& 0.5			\\
$g$								& maximum number of nearest neighbor customers evaluated by local search  	& 40			\\
$I_{0}$						& allowed HGS iterations without improvement for the first sub-problem      & 10000		\\
$I_{n}$						& allowed HGS iterations without improvement for subsequent sub-problems		& 500			\\ \hline
$n_{max}$					& maximum number of sub-problems solved by HRS   														& 50		\\ \hline
\end{tabular}%
}
\caption{Parameter settings of IHGS}
\label{params}
\end{table}

\subsection{Performance Metrics}
\label{4-2}
Rather than a single best known solution, MO algorithms must be judged on their approximation of a reference set, ideally the set of Pareto-optimal solutions to each instance. In order to generate a large and robust reference set that is as exhaustive as possible, we embedded our adapted HGS into a classical $\epsilon$-constraint framework, set the $\epsilon$ parameter to $0.01$, and solved every instance 10 times with $I_{n}$ quadrupled to $2000$. The final best known reference sets were formed by the non-dominated union of all solutions found in all our test runs, all solutions reported by \cite{oyola2014} and \cite{lacomme2015}, as well as the cost-optimal solutions reported in the literature.

It is generally not trivial how to assess whether or to what extent the Pareto set obtained with one algorithm is ``better'' than the set obtained by another algorithm. A comprehensive discussion about quality indicators is beyond the scope of this article, so the reader is referred to \cite{zitzler2003} for a detailed and rigorous analysis of the subject. We will evaluate solution quality based on the commonly used hypervolume and multiplicative unary epsilon indicators.

\paragraph{Hypervolume}
The hypervolume measures the volume of the objective space dominated by an approximation set. It is strictly monotonic (i.e.~sensitive to all additional non-dominated solutions), and it reaches its maximal value only when all Pareto-optimal solutions have been found -- it is thus far the only known MO quality measure which satisfies both of these properties \citep{zitzler2007}. 

\paragraph{Multiplicative Unary Epsilon}
This indicator measures the minimal value $\epsilon$ by which the objective values of all solutions in a reference set (e.g.~Pareto-optimal or best known) would have to be multiplied (minimization objectives) or divided (maximization objectives) so that all solutions in the reference set are dominated by a given approximation set. Hence, lower values are better, and if the approximation set consists of all reference solutions, then the optimal indicator value is exactly 1.

\bigskip
Loosely speaking, the hypervolume provides an average-case perspective, while the unary epsilon indicator, based solely on the largest gap to the reference set, gives a better worst-case perspective. Using both metrics, we aim to provide a more nuanced view. We also report solution set cardinality and CPU time.




The reference point for the hypervolume is set to the maximum encountered objective values, incremented by 1\% of the objective function ranges so that extreme solutions also contribute to the measure. For ease of interpretation, we report the hypervolume as a percentage of the reference hypervolume. Since the unary epsilon indicator is sensitive to scaling, we normalize the objective values of all solutions to the range [1,2] in order to normalize the unary epsilon to a fixed range and hence make performance on different instances more comparable. 

IHGS was coded in C++ and all computational experiments were conducted on a desktop computer using a single thread of an Intel Core i7-4790~3.60~GHz processor. The best known reference sets, reference points, as well as detailed computational results, are included with the online version of this article.

\subsection{Comparison with the Classical $\epsilon$-Constraint Framework}
\label{5-1}

The HRS framework generalizes the classical ECM, and so exactly the same solver can be embedded in both algorithms. Based on the theoretical discussions in Section \ref{3-0}, we are interested in evaluating also empirically the added value of HRS compared to the classical ECM, with respect to overall approximation quality, as well as convergence and anytime behavior.

The classical ECM requires the specification of the $\epsilon$ parameter. As discussed in Section \ref{3-0}, the ``best'' setting is instance-dependent. Considering that the classical ECM is meant to be exhaustive, we simply set $\epsilon = 0.01$ for all instances, in order to provide a fair basis of comparison. 


We solve the 20 E instances used by \cite{oyola2014} and the 14 CMT instances used by \cite{lacomme2015}. Each instance is solved 10 times with both versions of IHGS and the same parameter settings.

\begin{figure} [t]
	\centering
	\vspace{-0.6cm}
		\includegraphics[width=0.94\textwidth]{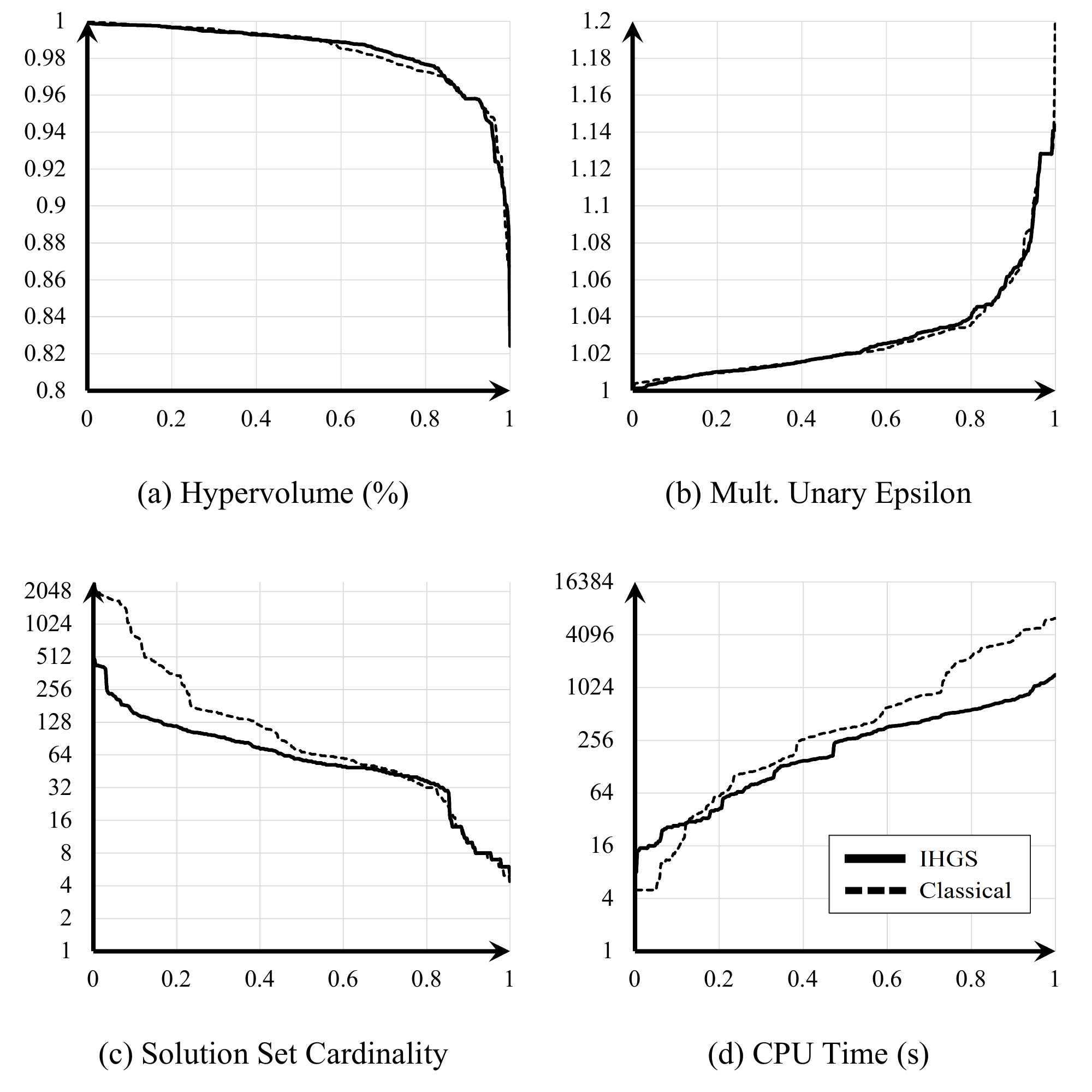}
	\caption{Attained final approximation set quality ($y$-axis) by a given percentage of all test runs ($x$-axis) ($n=340$).}
	\label{achievement}
\end{figure}

\paragraph{Final Approximation Quality}
In order to succinctly summarize the performance of both algorithm variants on all examined instances, we present in Figures \ref{achievement}a to \ref{achievement}d the achievement function for each performance metric: the $x$-axis indicates the percentage of test runs which achieve a performance \textit{at least as good} as the corresponding $y$-axis indicator value. Each data point corresponds to a single test run, so there is no averaging effect and all worst and best case scenarios are represented.

In terms of final approximation quality, we observe that both approaches attain hypervolume and unary epsilon values very close to the best known references -- 95\% of all test runs achieve a hypervolume of at least 95\% and unary epsilon values of no more than $1.1$, usually much lower. This shows that the classical ECM -- despite its drawbacks -- is already a good framework to build upon even with heuristics. Although the classical ECM variant identifies significantly more solutions on certain instances, in absolute terms the HRS framework still generates well over 100 solutions in those cases, which is arguably sufficient for an approximation, especially in light of the other quality indicators. We can therefore conclude that performing an exhaustive search with the classical ECM offers only very marginal benefits in terms of the quality of the generated Pareto set approximation.

When it comes to computational effort however, we observe that the more exhaustive exploration strategy of the classical ECM comes at the cost of a significantly higher computational burden. On all but the smallest instances, the classical ECM requires between two to four times as much computational effort (note the logarithmic scale of Figure \ref{achievement}-d).

\newpage

We emphasize again that although a larger $\epsilon$ value would speed up the search, there is no single value that works well for all instances. On these benchmarks, the values of the equity objective have ranges as small as $17$ to over $11000$, the number of solutions in the reference sets ranges from less than $10$ to over $5000$, and the distribution of these solutions in the objective space is highly erratic, to the point that sometimes half the non-dominated solutions are contained within only a fraction of the feasible range of the equity objective. This variability precludes a single ``best'' $\epsilon$ value ex ante.

Overall, these results demonstrate that the HRS framework generates solution sets whose final quality is virtually the same as those obtained with the most exhaustive version of the classical ECM, within only a fraction of the time. However, it is unclear to what extent an earlier termination of either method would affect these observations. Hence we focus our further attention on the differences in convergence.

\paragraph{Convergence Behavior}
In order to assess the added value of HRS on anytime behavior as well as the impact of using a population-based single-objective solver (\textit{solution harvesting}, cf.~Section \ref{3-3}), the current approximation set (with and without harvested solutions) was recorded after each $\epsilon$-constraint sub-problem. This allows for an analysis of performance as a function of CPU effort.

\begin{landscape}
\begin{figure} [!t]
	\centering
	\vspace{-0.7cm}
		\includegraphics[width=1.50\textwidth]{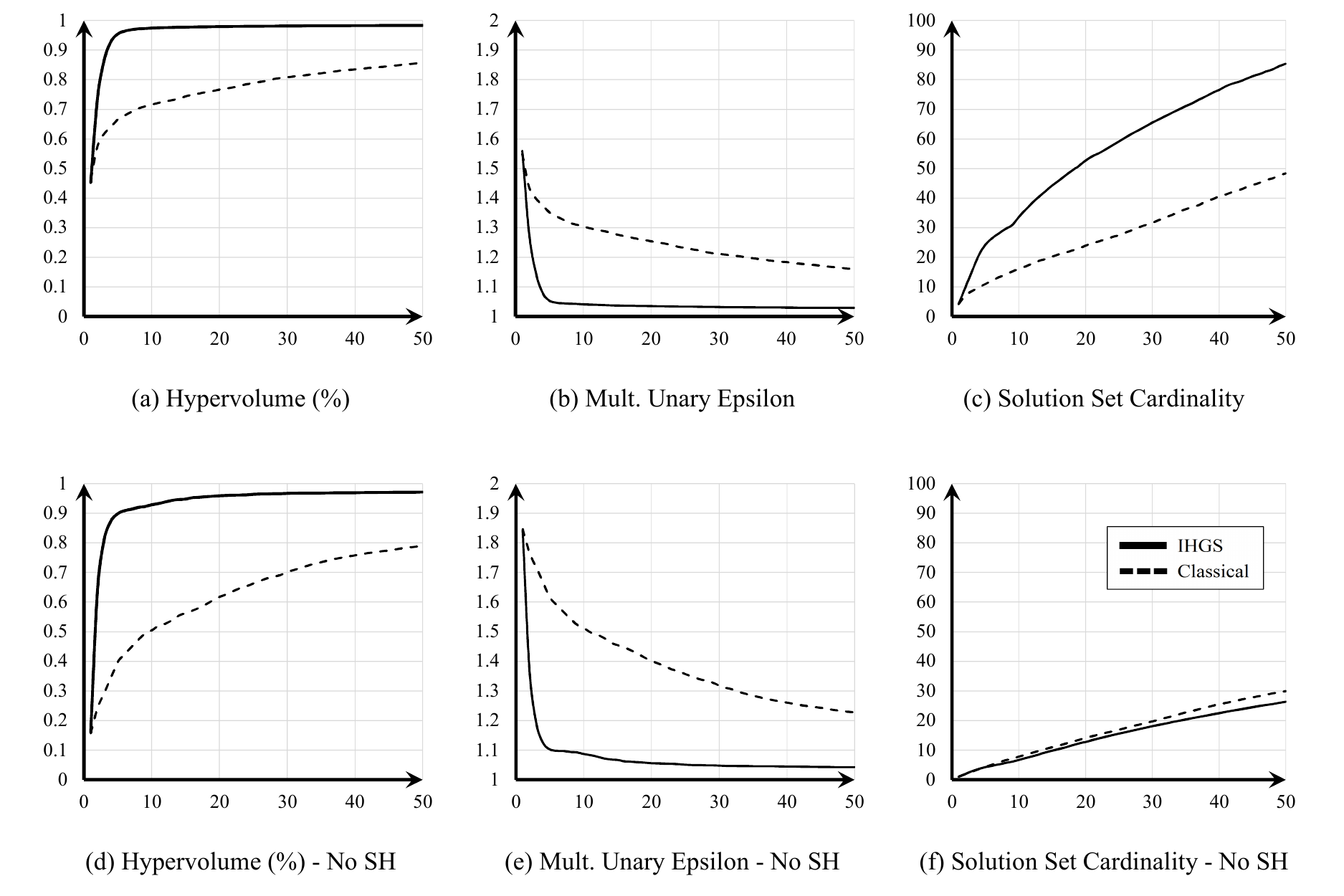}
	\caption{Convergence of performance indicator values ($y$-axis) as a function of the number of sub-problems solved ($x$-axis). SH stands for solution harvesting, i.e.~returning a population of potentially non-dominated solutions instead of a single cost-minimum per sub-problem.}
	\label{anytime+harvesting}
\end{figure}
\end{landscape}

Figures \ref{anytime+harvesting}a to \ref{anytime+harvesting}c present for both algorithm versions the development of the hypervolume, unary epsilon, and cardinality metrics as a function of the number of sub-problems solved, averaged over all instances and test runs. Figures \ref{anytime+harvesting}d to \ref{anytime+harvesting}f visualize the same comparison without solution harvesting, i.e.~the solver returns only the single cost-minimizing solution to each $\epsilon$-constraint sub-problem.

From the data we can clearly see that the binary search of the HRS framework converges significantly faster than the classical ECM, with respect to all three solution quality metrics. In fact, we find that setting $n_{max}=50$ was rather conservative, as on average the hypervolume and unary epsilon converge nearly to their final values after solving only 10 sub-problems. In contrast, the classical ECM is still far from producing a reasonable approximation at that point in time. This casts the previously discussed differences in computational effort in a more dramatic light, as HRS could have been terminated much earlier with marginal effect on approximation quality. The same cannot be said of the classical ECM.

With respect to the impact of solution harvesting, we make two observations: First, it speeds up the convergence of the classical ECM much more than that of the HRS framework. Since consecutive solutions found with the classical approach are generally very close in the objective space, the impact of more diverse solutions is large. The opposite is true for HRS. Second, solution harvesting has little effect on \textit{final} approximation quality in terms of the attained hypervolume and unary epsilon values, despite greater cardinality. This reveals that harvested solutions tend to be similar to the respective sub-problem cost-optima, and hence do not contribute noticeably to approximating the trade-off structure. It also has methodological implications, as it demonstrates that single-solution-based solvers can be used to generate high-quality Pareto set approximations. However, there is no significant downside to using solution harvesting, and single-solution-based methods can also employ the concept, e.g.~by storing encountered local optima.


\bigskip
Based on the above analyses, we have demonstrated that HRS improves upon the classical ECM by converging much more rapidly while still achieving virtually the same high approximation quality as an exhaustive search strategy. In the next two sections, we contrast the performance of HRS/IHGS with the current state-of-the-art multi-objective metaheuristics.

\subsection{Comparison with GRASP-ASP}
\label{5-2}

\begin{figure} [!b]
	\centering
		\includegraphics[width=1.0\textwidth]{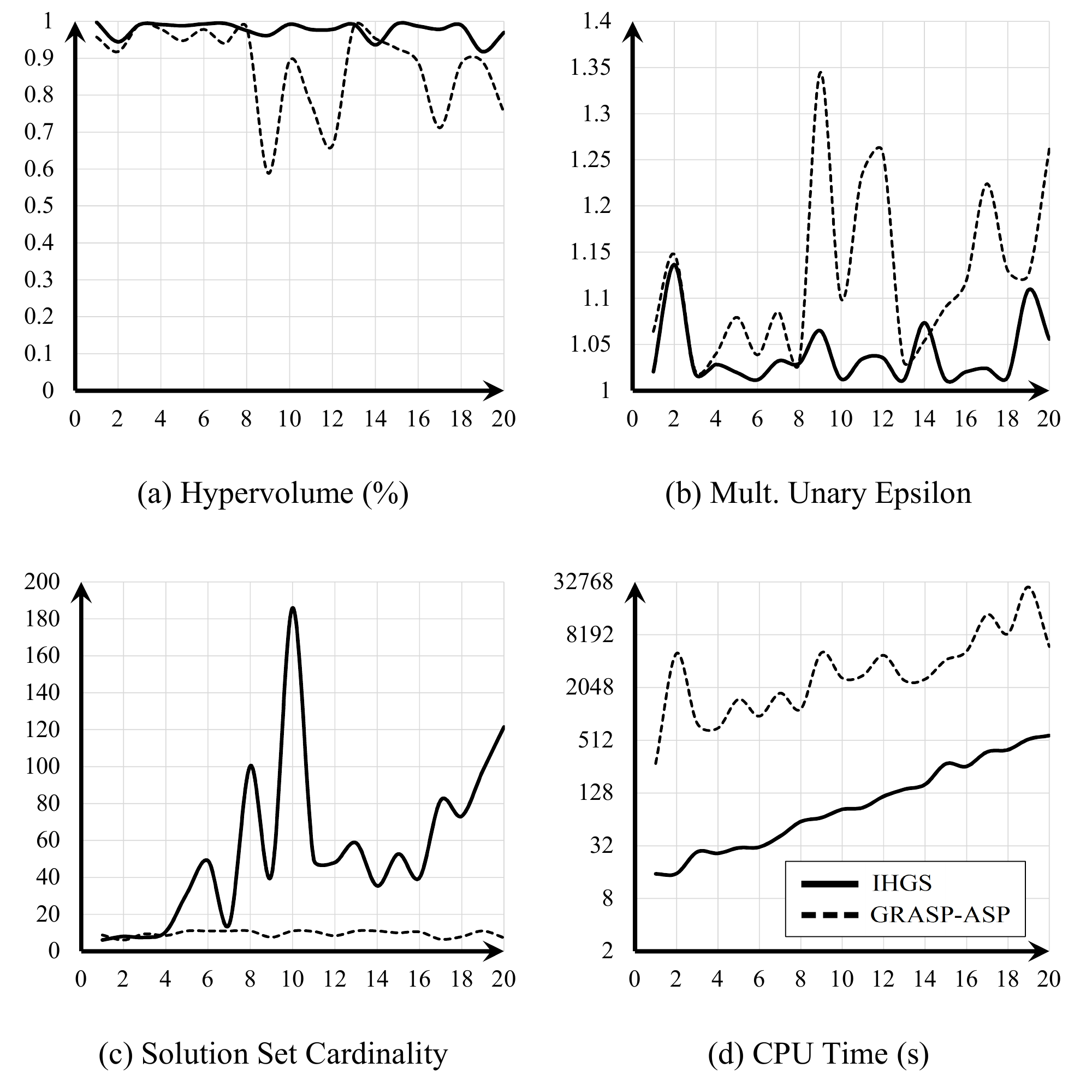}
	\caption{Average performance per instance of IHGS vs.~GRASP-ASP}
	\label{E_visual}
\end{figure}

We compare the performance of IHGS with that of GRASP-ASP on the 20 E instances used by \cite{oyola2014}, following the same conventions: the fleet size is fixed to that of the cost-optimal solution, each unused vehicle is counted in the equity objective with a workload of 0, and only 2-optimal solutions are considered for inclusion in the archive. Tables \ref{e-ihgs} and \ref{e-grasp-asp} report for each instance and algorithm the attained hypervolume, unary epsilon, Pareto set cardinality, and CPU time. The tables list for each instance and performance metric the maximum, average, and minimum values obtained out of 10 test runs. Figures \ref{E_visual}a to \ref{E_visual}d visualize average performance per instance.




{

\begin{table}[!h]
\centering
\resizebox{\textwidth}{!}{%
\begin{tabular}{C|CCC|CCC|CCC|CCC}
\hline
\textbf{} & \multicolumn{3}{c|}{\textbf{Hypervolume (\%)}} & \multicolumn{3}{c|}{\textbf{Unary Epsilon}} & \multicolumn{3}{c|}{\textbf{Cardinality}} & \multicolumn{3}{c}{\textbf{CPU Time (s)}} \\
\textbf{Inst.}	& \textit{Max} & \textit{Ave} & \textit{Min} & \textit{Max} & \textit{Ave} & \textit{Min} & \textit{Max} & \textit{Ave} & \textit{Min} & \textit{Max} & \textit{Ave} & \textit{Min} \\ \hline
\textbf{1} & 99.9 & 99.8 & 99.1 & 1.020 & 1.020 & 1.020 & 7 & 6.1 & 5 & 18 & 15.3 & 8 \\
\textbf{2} & 95.8 & 94.5 & 82.4 & 1.211 & 1.137 & 1.128 & 8 & 8.0 & 8 & 16 & 15.4 & 14 \\
\textbf{3} & 99.5 & 99.1 & 98.8 & 1.021 & 1.019 & 1.014 & 9 & 7.5 & 7 & 30 & 27.2 & 25 \\
\textbf{4} & 99.8 & 99.2 & 98.1 & 1.041 & 1.028 & 1.021 & 13 & 10.6 & 9 & 28 & 26.2 & 24 \\
\textbf{5} & 99.4 & 98.8 & 97.7 & 1.033 & 1.020 & 1.016 & 34 & 31.3 & 28 & 34 & 30.3 & 20 \\
\textbf{6} & 99.5 & 99.3 & 98.9 & 1.017 & 1.012 & 1.010 & 49 & 49.0 & 49 & 33 & 30.9 & 29 \\
\textbf{7} & 99.8 & 99.5 & 99.4 & 1.034 & 1.032 & 1.023 & 17 & 14.5 & 14 & 43 & 41.0 & 39 \\
\textbf{8} & 98.7 & 97.5 & 95.5 & 1.052 & 1.030 & 1.017 & 104 & 100.4 & 93 & 66 & 60.3 & 54 \\
\textbf{9} & 98.1 & 96.2 & 94.6 & 1.084 & 1.065 & 1.046 & 43 & 41.1 & 39 & 75 & 66.9 & 61 \\
\textbf{10} & 99.4 & 99.2 & 99.1 & 1.014 & 1.013 & 1.011 & 207 & 186.0 & 171 & 93 & 83.2 & 72 \\
\textbf{11} & 99.3 & 97.8 & 96.1 & 1.052 & 1.034 & 1.015 & 53 & 49.7 & 47 & 96 & 87.2 & 80 \\
\textbf{12} & 98.6 & 97.8 & 97.2 & 1.046 & 1.036 & 1.028 & 54 & 48.1 & 41 & 141 & 117.2 & 66 \\
\textbf{13} & 99.6 & 99.2 & 98.8 & 1.014 & 1.011 & 1.006 & 63 & 58.7 & 54 & 151 & 140.9 & 132 \\
\textbf{14} & 96.6 & 93.7 & 89.6 & 1.116 & 1.074 & 1.041 & 38 & 35.4 & 33 & 168 & 160.2 & 145 \\
\textbf{15} & 99.8 & 99.3 & 98.8 & 1.019 & 1.012 & 1.007 & 57 & 52.7 & 48 & 302 & 275.3 & 260 \\
\textbf{16} & 99.5 & 98.7 & 98.0 & 1.027 & 1.020 & 1.014 & 43 & 39.8 & 37 & 302 & 257.3 & 228 \\
\textbf{17} & 98.5 & 97.9 & 97.5 & 1.028 & 1.024 & 1.017 & 90 & 81.7 & 68 & 465 & 376.7 & 252 \\
\textbf{18} & 99.6 & 99.0 & 98.2 & 1.022 & 1.015 & 1.008 & 83 & 73.2 & 66 & 431 & 398.3 & 373 \\
\textbf{19} & 95.3 & 91.8 & 88.3 & 1.151 & 1.109 & 1.062 & 106 & 97.9 & 90 & 569 & 525.8 & 493 \\
\textbf{20} & 97.8 & 97.0 & 95.8 & 1.070 & 1.056 & 1.046 & 136 & 121.5 & 110 & 621 & 578.6 & 524 \\	\hline
\textbf{Ave} & \textbf{98.7} & \textbf{97.8} & \textbf{96.1} & \textbf{1.054} & \textbf{1.038} & \textbf{1.028} & \textbf{60.7} & \textbf{55.7} & \textbf{50.9} & \textbf{184.1} & \textbf{165.7} & \textbf{144.9}	\\ \hline
\end{tabular}%
}
\caption{Performance of IHGS on the E instances}
\label{e-ihgs}
\end{table}
}

{
\newcolumntype{C}{>{\centering\arraybackslash}p{2.5em}}

\begin{table}[!t]
\centering
\resizebox{\textwidth}{!}{%
\begin{tabular}{C|CCC|CCC|CCC|CCC}
\hline
\textbf{} & \multicolumn{3}{c|}{\textbf{Hypervolume (\%)}} & \multicolumn{3}{c|}{\textbf{Unary Epsilon}} & \multicolumn{3}{c|}{\textbf{Cardinality}} & \multicolumn{3}{c}{\textbf{CPU Time (s)}} \\
\textbf{Inst.}	& \textit{Max} & \textit{Ave} & \textit{Min} & \textit{Max} & \textit{Ave} & \textit{Min} & \textit{Max} & \textit{Ave} & \textit{Min} & \textit{Max} & \textit{Ave} & \textit{Min} \\ \hline
\textbf{1} & 99.8 & 95.8 & 73.5 & 1.277 & 1.064 & 1.013 & 11 & 8.8 & 7 & - & 277 & - \\
\textbf{2} & 100.0 & 91.8 & 80.9 & 1.236 & 1.148 & 1.000 & 8 & 6.1 & 4 & - & 5005 & - \\
\textbf{3} & 99.7 & 99.1 & 98.5 & 1.031 & 1.022 & 1.014 & 11 & 9.4 & 7 & - & 799 & - \\
\textbf{4} & 98.9 & 97.9 & 94.9 & 1.058 & 1.040 & 1.035 & 9 & 8.5 & 8 & - & 702 & - \\
\textbf{5} & 97.8 & 94.8 & 87.7 & 1.152 & 1.080 & 1.029 & 11 & 11 & 11 & - & 1503 & - \\
\textbf{6} & 98.3 & 97.8 & 97.4 & 1.049 & 1.039 & 1.026 & 11 & 11 & 11 & - & 963 & - \\
\textbf{7} & 99.2 & 94.0 & 90.1 & 1.102 & 1.085 & 1.034 & 11 & 11 & 11 & - & 1765 & - \\
\textbf{8} & 98.5 & 98.3 & 98.1 & 1.037 & 1.031 & 1.024 & 11 & 11 & 11 & - & 1161 & - \\
\textbf{9} & 77.7 & 58.9 & 29.8 & 1.500 & 1.345 & 1.217 & 11 & 7.6 & 3 & - & 5111 & - \\
\textbf{10} & 95.5 & 89.3 & 67.9 & 1.266 & 1.100 & 1.038 & 11 & 11 & 11 & - & 2632 & - \\
\textbf{11} & 84.2 & 77.7 & 70.2 & 1.288 & 1.233 & 1.175 & 11 & 10.7 & 8 & - & 2804 & - \\
\textbf{12} & 78.9 & 66.4 & 46.8 & 1.481 & 1.258 & 1.148 & 11 & 8.4 & 3 & - & 4771 & - \\
\textbf{13} & 99.3 & 98.6 & 98.0 & 1.043 & 1.032 & 1.023 & 11 & 11 & 11 & - & 2465 & - \\
\textbf{14} & 98.2 & 95.4 & 92.8 & 1.083 & 1.054 & 1.036 & 11 & 11 & 11 & - & 2542 & - \\
\textbf{15} & 96.1 & 92.7 & 88.4 & 1.128 & 1.090 & 1.044 & 11 & 10 & 7 & - & 4203 & - \\
\textbf{16} & 92.6 & 88.8 & 85.5 & 1.154 & 1.118 & 1.095 & 11 & 10.5 & 9 & - & 5290 & - \\
\textbf{17} & 76.9 & 71.2 & 67.0 & 1.268 & 1.224 & 1.173 & 10 & 6.5 & 4 & - & 13917 & - \\
\textbf{18} & 91.0 & 88.6 & 85.6 & 1.149 & 1.130 & 1.112 & 10 & 7.9 & 6 & - & 8401 & - \\
\textbf{19} & 91.5 & 89.1 & 81.6 & 1.174 & 1.126 & 1.103 & 11 & 11 & 11 & - & 28721 & - \\
\textbf{20} & 81.7 & 75.3 & 64.6 & 1.341 & 1.262 & 1.215 & 11 & 7.3 & 5 & - & 5998 & - \\	\hline
\textbf{Ave} & \textbf{92.8} & \textbf{88.1} & \textbf{80.0} & \textbf{1.191} & \textbf{1.124} & \textbf{1.078} & \textbf{10.7} & \textbf{9.5} & \textbf{8.0} & \textbf{-} & \textbf{4951} & \textbf{-}	\\ \hline
\end{tabular}%
}
\caption{Performance of GRASP-ASP on the E instances}
\label{e-grasp-asp}
\end{table}
}

\paragraph{Hypervolume}
On the smallest instances, there is no significant difference between IHGS and GRASP-ASP, but as the number of customers and vehicles increases, the gap between the two algorithms widens noticeably. Considering all 20 instances together, IHGS achieves an average hypervolume of over 97\% compared to 88\% for GRASP-ASP. In addition, IHGS is significantly more robust: over all instances, the average \textit{worst} performance of IHGS -- 96.1\% -- is higher than the average \textit{best} performance of GRASP-ASP, at 92.8\%.

\paragraph{Unary Epsilon}
IHGS generally outperforms GRASP-ASP with respect to the epsilon metric, and the gap between the two algorithms also increases with instance size. The average indicator value over all instances is only 1.038 for IHGS, compared to 1.124 for GRASP-ASP. As with the hypervolume, the average worst unary epsilon value of IHGS -- only 1.054 -- is better than the average best attained by GRASP-ASP, at 1.078.

\paragraph{Cardinality}
Overall, IHGS generates around five times as many non-dominated solutions as GRASP-ASP. On its own, a larger cardinality does not imply a better approximation, but when interpreted in conjunction with the unary epsilon indicator, it becomes clear that not only does IHGS find \textit{more} trade-off solutions, \textit{all} of them are closer to the reference set even in the worst case. This demonstrates further the robustness of IHGS.

\paragraph{CPU Time}
The CPU times of different algorithms can vary widely due to factors beyond algorithm design, such as the programming language, the processor, and so on. However, it is unlikely that these factors alone could account for a difference of several orders of magnitude: IHGS is around 16 times faster than GRASP-ASP on all examined instances.

\bigskip
Overall, we can conclude that IHGS generates Pareto set approximations of markedly higher quality and with many more trade-off solutions, and achieves this with only a small fraction of the computational effort used by GRASP-ASP. Looking at the benchmark set as a whole, and at both average and worst case performance, it is fair to say that IHGS outperforms GRASP-ASP with respect to all four examined performance metrics.


\subsection{Comparison with MSSPR}
\label{5-3}

\begin{figure} [!t]
	\centering
		\includegraphics[width=1.0\textwidth]{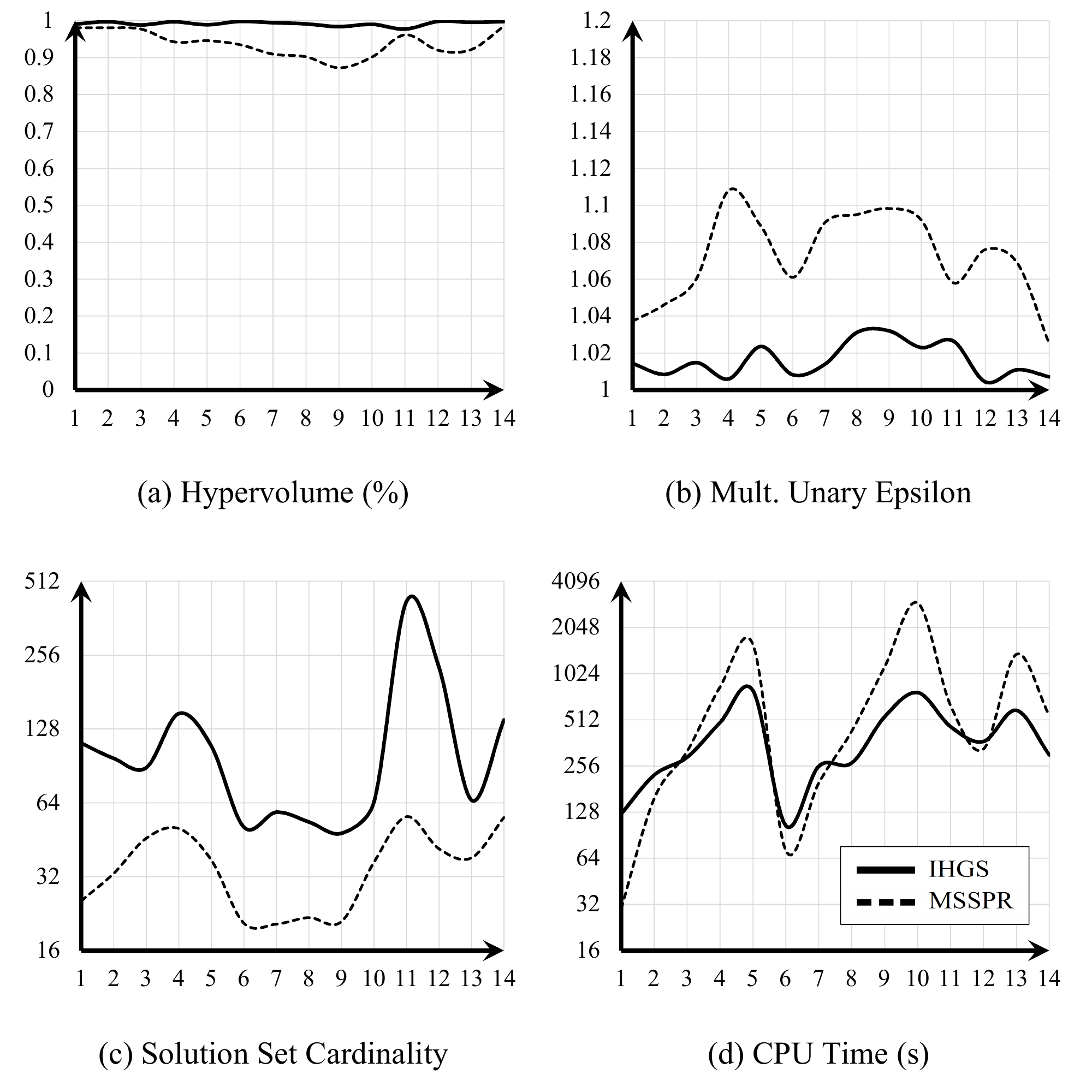}
	\caption{Average performance per instance of IHGS vs.~MSSPR}
	\label{CMT_visual}
\end{figure}

Our comparison with the MSSPR of \cite{lacomme2015} is based on the 14 CMT instances introduced by \cite{CMT}. As in \cite{lacomme2015}, neither a fixed fleet size nor the 2-optimality constraint are imposed. However, we point out that the fleet size convention appears to be relevant only for instances CMT5 and CMT14, in which sometimes the equity objective can be improved by adding or removing one vehicle. As before, Tables \ref{c-ihgs} and \ref{c-msspr} report for both algorithms the maximum, average, and minimum values (per instance) for all examined performance metrics, and the average case performance is visualized in Figures \ref{CMT_visual}a to \ref{CMT_visual}d. The results for MSSPR are based on 5 runs per instance.



{
\newcolumntype{C}{>{\centering\arraybackslash}p{2.5em}}

\begin{table}[!t]
\centering
\resizebox{\textwidth}{!}{%
\begin{tabular}{C|CCC|CCC|CCC|CCC}
\hline
\textbf{} & \multicolumn{3}{c|}{\textbf{Hypervolume (\%)}} & \multicolumn{3}{c|}{\textbf{Unary Epsilon}} & \multicolumn{3}{c|}{\textbf{Cardinality}} & \multicolumn{3}{c}{\textbf{CPU Time (s)}} \\
\textbf{Inst.}	& \textit{Max} & \textit{Ave} & \textit{Min} & \textit{Max} & \textit{Ave} & \textit{Min} & \textit{Max} & \textit{Ave} & \textit{Min} & \textit{Max} & \textit{Ave} & \textit{Min} \\ \hline
\textbf{1} & 99.9 & 99.1 & 97.3 & 1.029 & 1.014 & 1.002 & 133 & 112.2 & 93 & 171 & 154.0 & 124 \\
\textbf{2} & 99.8 & 99.7 & 99.7 & 1.011 & 1.008 & 1.007 & 120 & 97.1 & 64 & 369 & 333.8 & 267 \\
\textbf{3} & 99.9 & 98.9 & 96.4 & 1.037 & 1.015 & 1.005 & 121 & 89.0 & 65 & 463 & 413.7 & 296 \\
\textbf{4} & 99.8 & 99.7 & 99.4 & 1.008 & 1.006 & 1.005 & 164 & 148.0 & 127 & 925 & 803.7 & 640 \\
\textbf{5} & 99.3 & 98.9 & 98.0 & 1.037 & 1.024 & 1.019 & 151 & 109.1 & 71 & 1200 & 1101.3 & 939 \\
\textbf{6} & 99.9 & 99.8 & 99.8 & 1.014 & 1.008 & 1.005 & 63 & 51.0 & 42 & 157 & 145.9 & 130 \\
\textbf{7} & 99.7 & 99.5 & 99.3 & 1.018 & 1.014 & 1.010 & 77 & 58.6 & 46 & 322 & 294.6 & 243 \\
\textbf{8} & 99.4 & 99.1 & 98.8 & 1.036 & 1.031 & 1.028 & 75 & 53.5 & 43 & 427 & 379.7 & 337 \\
\textbf{9} & 98.8 & 98.4 & 97.7 & 1.040 & 1.032 & 1.028 & 66 & 48.1 & 37 & 831 & 738.7 & 630 \\
\textbf{10} & 99.4 & 99.0 & 98.6 & 1.029 & 1.023 & 1.019 & 83 & 64.4 & 44 & 1430 & 1254.0 & 993 \\
\textbf{11} & 99.6 & 97.7 & 94.5 & 1.059 & 1.027 & 1.009 & 488 & 423.5 & 391 & 730 & 635.0 & 527 \\
\textbf{12} & 99.9 & 99.8 & 99.8 & 1.005 & 1.004 & 1.004 & 258 & 228.7 & 206 & 560 & 481.6 & 409 \\
\textbf{13} & 99.9 & 99.6 & 99.4 & 1.015 & 1.011 & 1.009 & 94 & 65.7 & 41 & 809 & 713.9 & 590 \\
\textbf{14} & 99.9 & 99.8 & 99.7 & 1.008 & 1.007 & 1.006 & 160 & 139.4 & 130 & 583 & 548.2 & 440 \\	\hline
\textbf{Ave} & \textbf{99.6} & \textbf{99.2} & \textbf{98.4} & \textbf{1.025} & \textbf{1.016} & \textbf{1.011} & \textbf{146.6} & \textbf{120.6} & \textbf{100.0} & \textbf{641.2} & \textbf{571.3} & \textbf{468.9}	\\ \hline
\end{tabular}%
}
\caption{Performance of IHGS on the CMT instances}
\label{c-ihgs}
\end{table}
}

{
\newcolumntype{C}{>{\centering\arraybackslash}p{2.5em}}

\begin{table}[!t]
\centering
\resizebox{\textwidth}{!}{%
\begin{tabular}{C|CCC|CCC|CCC|CCC}
\hline
\textbf{} & \multicolumn{3}{c|}{\textbf{Hypervolume (\%)}} & \multicolumn{3}{c|}{\textbf{Unary Epsilon}} & \multicolumn{3}{c|}{\textbf{Cardinality}} & \multicolumn{3}{c}{\textbf{CPU Time (s)}} \\
\textbf{Inst.}	& \textit{Max} & \textit{Ave} & \textit{Min} & \textit{Max} & \textit{Ave} & \textit{Min} & \textit{Max} & \textit{Ave} & \textit{Min} & \textit{Max} & \textit{Ave} & \textit{Min} \\ \hline
\textbf{1} & 99.0 & 98.1 & 96.8 & 1.059 & 1.037 & 1.020 & 35 & 25.5 & 17 & - & 30 & - \\
\textbf{2} & 99.1 & 98.1 & 96.8 & 1.059 & 1.046 & 1.036 & 47 & 33.0 & 17 & - & 156 & - \\
\textbf{3} & 99.1 & 97.7 & 96.1 & 1.089 & 1.060 & 1.036 & 60 & 45.9 & 34 & - & 318 & - \\
\textbf{4} & 97.2 & 94.3 & 90.7 & 1.147 & 1.108 & 1.070 & 60 & 50.3 & 40 & - & 834 & - \\
\textbf{5} & 97.6 & 94.6 & 90.7 & 1.147 & 1.089 & 1.037 & 60 & 37.4 & 23 & - & 1590 & - \\
\textbf{6} & 97.6 & 93.5 & 84.6 & 1.115 & 1.061 & 1.037 & 27 & 20.7 & 12 & - & 72 & - \\
\textbf{7} & 95.1 & 90.9 & 84.6 & 1.115 & 1.091 & 1.052 & 28 & 20.5 & 12 & - & 198 & - \\
\textbf{8} & 93.6 & 90.2 & 85.7 & 1.112 & 1.095 & 1.079 & 32 & 21.8 & 12 & - & 432 & - \\
\textbf{9} & 90.2 & 87.3 & 83.5 & 1.127 & 1.098 & 1.079 & 32 & 21.0 & 12 & - & 1134 & - \\
\textbf{10} & 95.5 & 90.2 & 83.5 & 1.127 & 1.092 & 1.056 & 74 & 36.6 & 19 & - & 2964 & - \\
\textbf{11} & 98.4 & 96.2 & 93.4 & 1.099 & 1.058 & 1.036 & 75 & 56.2 & 34 & - & 636 & - \\
\textbf{12} & 98.4 & 92.0 & 82.4 & 1.147 & 1.076 & 1.036 & 75 & 41.6 & 13 & - & 330 & - \\
\textbf{13} & 99.0 & 92.2 & 82.4 & 1.147 & 1.069 & 1.023 & 63 & 38.2 & 13 & - & 1368 & - \\
\textbf{14} & 99.0 & 98.6 & 98.1 & 1.029 & 1.025 & 1.023 & 63 & 55.8 & 43 & - & 546 & - \\	\hline
\textbf{Ave} & \textbf{95.0} & \textbf{93.5} & \textbf{92.0} & \textbf{1.109} & \textbf{1.072} & \textbf{1.044} & \textbf{52.2} & \textbf{36.0} & \textbf{21.5} & \textbf{-} & \textbf{757.7} & \textbf{-}	\\	\hline
\end{tabular}%
}
\caption{Performance of MSSPR on the CMT instances}
\label{c-msspr}
\end{table}
}

\paragraph{Hypervolume}
IHGS achieves an average hypervolume per instance that is consistently higher than that achieved by MSSPR. In fact, for nearly every instance the worst performance obtained with IHGS is better than the corresponding best performance of MSSPR. IHGS is also more robust: the lowest hypervolume attained by IHGS in all test runs is around 94\%, in contrast to 82\% for MSSPR.

\paragraph{Unary Epsilon}
Although the unary epsilon values are low for both algorithms, the average values obtained with IHGS are always lower than those of MSSPR. As with the hypervolume, the worst unary epsilon value obtained with IHGS is often better than the best value attained by MSSPR on the same instance.

\paragraph{Cardinality}
IHGS finds significantly larger approximation sets on every CMT instance: on average IHGS identifies over three times as many trade-off solutions as MSSPR. As pointed out in the previous section, finding more non-dominated solutions is not an end in itself. However, the hypervolume and unary epsilon metrics indicate that the Pareto sets generated with IHGS are not just larger, they are also of higher quality. 

\paragraph{CPU Time}
The computational effort required by IHGS and MSSPR appears to be largely similar. However, the CPU times of MSSPR seem to grow noticeably faster with increasing instance size. In addition, MSSPR appears to be more sensitive to additional constraints, as its computational burden increases when the duration constraint is imposed (CMT6 to CMT10 are identical to the previous five, except for the duration constraint). Due to the many factors that affect CPU time, it is difficult to make any further conclusions, but based on the convergence behavior analyzed in Section \ref{5-1}, we remark that IHGS could have been terminated after half the time and with only marginal impact on final approximation quality.

\bigskip
Overall, IHGS clearly outperforms MSSPR in terms of solution set quality with computational effort that is of a comparable order of magnitude.

\begin{figure} [!b]
	\centering
		\includegraphics[width=1.0\textwidth]{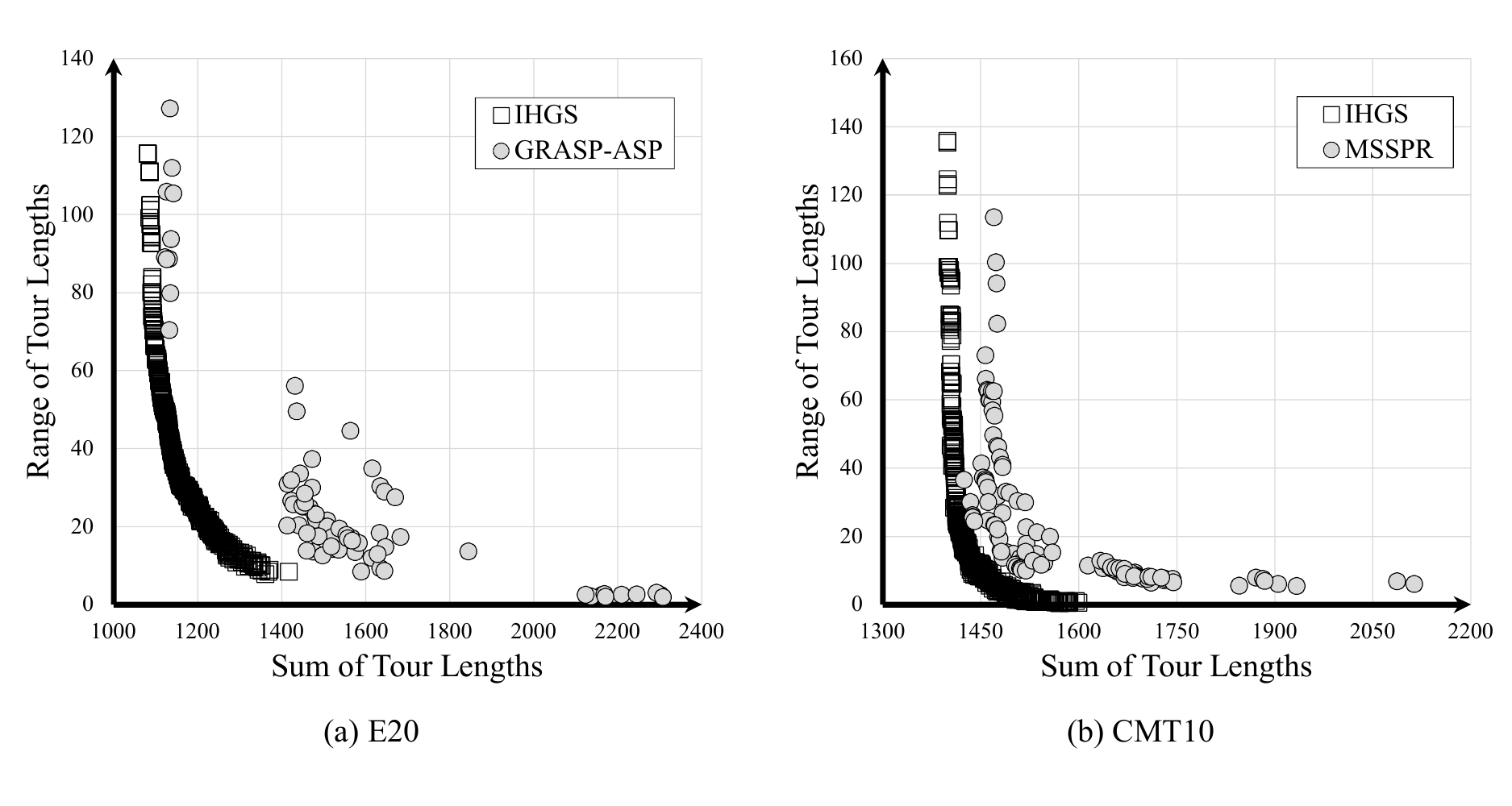}
	\caption{All generated Pareto sets on two representative instances}
	\label{fronts}
\end{figure}

\bigskip
We close our empirical study with a comparison of the actual Pareto sets generated by the examined algorithms on the two largest instances: E20 with 100 customers and 14 vehicles, and CMT10 with 199 customers, 18 vehicles, and a duration constraint. Figures \ref{fronts}a and \ref{fronts}b visualize the Pareto sets generated by IHGS, GRASP-ASP, and MSSPR on these instances. We plot the Pareto sets of all test runs in order to provide a more accurate impression of reliability.

From the figures we can see that (a) the quality of the solutions found with IHGS is more consistent, (b) the distribution of these solutions is more even, and that (c) IHGS solutions Pareto-dominate nearly all the solutions identified by the competing approaches. In fact, considering all instances and all pairwise comparisons of test runs, IHGS dominates on average around 75\% of all solutions found by GRASP-ASP and MSSPR. Based on the foregoing analyses, we believe it is fair to conclude that IHGS significantly outperforms the state-of-the-art approaches with respect to all examined performance metrics.

\section{Conclusion}
\label{6-conclusion}

In principle, the $\epsilon$-constraint method (ECM) offers a simple and direct way to leverage the many decades of research on single-objective heuristics to effectively solve multi-objective problems. Yet despite offering attractive theoretical properties and being straightforward to implement, it has seen very limited use with heuristics due to a number of unresolved challenges. Based on these observations, we proposed an adaptation -- Heuristic Rectangle Splitting (HRS) -- that retains the essential properties of the standard method while mitigating the main difficulties arising in a heuristic context.

\paragraph{Empirical Results}
We have demonstrated the practical potential of the HRS framework on the VRP with Route Balancing (VRPRB) by embedding a state-of-the-art VRP heuristic -- the Hybrid Genetic Search of \cite{vidal2012} -- as the underlying single-objective solver. Based on an extensive computational study considering various performance metrics, we have shown that HRS converges significantly faster than the classical ECM, and does so without compromising solution quality.

We then contrasted our approach with the current state-of-the-art multi-objective metaheuristics for the VRPRB: GRASP-ASP \citep{oyola2014} and MSSPR \citep{lacomme2015}. The results of this second empirical study demonstrate that HRS significantly outperforms both methods with respect to the hypervolume, multiplicative unary epsilon, and cardinality metrics, while requiring computational effort of a similar or lower order of magnitude. In addition, the HRS framework is problem-independent, simple to implement with only a CPU budget parameter, and usable with existing tailored heuristics. It also directly benefits from future progress on single-objective optimization.

\paragraph{Success Factors}
Although HRS is a generic framework, some heuristics will be better suited than others when used as the underlying solver. As a requirement, the heuristic must be flexible enough to handle either an additional constraint, or incorporate a penalty function. We believe the latter strategy offers more synergy with an $\epsilon$-constraint framework, since it allows solutions from previous steps to be more easily re-used after the $\epsilon$-constraint has been tightened. In addition, penalty functions are more generic and often improve the performance of heuristics by enabling the exploration of the infeasible solution space.

The other main success factor is the ability to store and exploit some form of search memory. This is critical for computational efficiency, as it allows the algorithm to warm-start from previous search states and more quickly converge to the next non-dominated solution. If multiple search states are kept in memory, then in later stages the algorithm can easily return to specific areas in the objective space to intensify the search where needed. Finally, our empirical analysis revealed that solution harvesting has only a marginal impact on overall approximation set quality, which suggests that population-based approaches are not necessarily superior to single-solution-based methods when embedded in an $\epsilon$-constraint framework such as HRS.

\paragraph{Future Research}
There are many promising avenues for further research, since the $\epsilon$-constraint framework has received limited attention in the literature on heuristic optimization. One issue that deserves more attention is the concept of minimal representative Pareto sets \citep{hamacher2007}: here the aim is to identify neither all nor as many as possible Pareto-optimal solutions, but rather a \textit{minimal} subset which provides a certain approximation quality according to one or more quality metrics. This is an area well-suited to heuristics in general, but also to the $\epsilon$-constraint framework in particular, since the latter can be used to select and precisely target specific regions in the objective space. Further methodological advances may also be made by considering the parallelization potential of $\epsilon$-constraint algorithms: since each sub-problem is entirely independent, a parallel search strategy combined with warm-starts could sharply reduce computational effort, which is a particularly important issue in multi-objective settings.


With respect to balanced VRPs, we remark that although the VRPRB is an informative and prototypical benchmark for purely numerical comparisons, recent work has revealed that it suffers from some modeling issues and should therefore be revisited before being used as the basis for richer problems \citep{matl2016}. Hence there is also room for introducing improved models and corresponding benchmarks. Last but not least, our empirical study revealed the strong performance of HRS in a VRP setting -- in light of these results, similar studies investigating the performance of HRS or related $\epsilon$-constraint-based heuristics on other problem classes appear highly promising.





\begin{spacing}{1}
\bibliographystyle{apa}
\bibliography{bibliography}

\begin{thebibliography}{}

\bibitem[\protect\astroncite{Boland et~al.}{2015}]{boland2015}
Boland, N., Charkhgard, H., and Savelsbergh, M. (2015).
\newblock A criterion space search algorithm for biobjective integer
  programming: The balanced box method.
\newblock {\em INFORMS Journal on Computing}, 27(4):735--754.

\bibitem[\protect\astroncite{Borgulya}{2008}]{borgulya2008}
Borgulya, I. (2008).
\newblock An algorithm for the capacitated vehicle routing problem with route
  balancing.
\newblock {\em Central European Journal of Operations Research},
  16(4):331--343.

\bibitem[\protect\astroncite{Boudia et~al.}{2007}]{boudia2007}
Boudia, M., Prins, C., and Reghioui, M. (2007).
\newblock An effective memetic algorithm with population management for the
  split delivery vehicle routing problem.
\newblock In {\em International Workshop on Hybrid Metaheuristics}, pages
  16--30. Springer.

\bibitem[\protect\astroncite{Christofides et~al.}{1979}]{CMT}
Christofides, N., Mingozzi, A., and Toth, P. (1979).
\newblock The vehicle routing problem.
\newblock In Christofides, N., Mingozzi, A., Toth, P., and Sandi, C., editors,
  {\em Combinatorial Optimization}, chapter~11, pages 315--338. John Wiley,
  Chichester.

\bibitem[\protect\astroncite{Hamacher et~al.}{2007}]{hamacher2007}
Hamacher, H.~W., Pedersen, C.~R., and Ruzika, S. (2007).
\newblock Finding representative systems for discrete bicriterion optimization
  problems.
\newblock {\em Operations Research Letters}, 35(3):336--344.

\bibitem[\protect\astroncite{Huang et~al.}{2004}]{huang2004}
Huang, B., Cheu, R.~L., and Liew, Y.~S. (2004).
\newblock {GIS} and genetic algorithms for hazmat route planning with security
  considerations.
\newblock {\em International Journal of Geographical Information Science},
  18(8):769--787.

\bibitem[\protect\astroncite{Jozefowiez et~al.}{2002}]{jozefowiez2002}
Jozefowiez, N., Semet, F., and Talbi, E.-G. (2002).
\newblock Parallel and hybrid models for multi-objective optimization:
  Application to the vehicle routing problem.
\newblock In {\em Parallel Problem Solving from Nature - PPSN VII}, pages
  271--280. Springer.

\bibitem[\protect\astroncite{Jozefowiez et~al.}{2006}]{jozefowiez2006}
Jozefowiez, N., Semet, F., and Talbi, E.-G. (2006).
\newblock Enhancements of {NSGA II} and its application to the vehicle routing
  problem with route balancing.
\newblock In {\em Artificial Evolution}, pages 131--142. Springer.

\bibitem[\protect\astroncite{Jozefowiez et~al.}{2007}]{jozefowiez2007}
Jozefowiez, N., Semet, F., and Talbi, E.-G. (2007).
\newblock Target aiming {P}areto search and its application to the vehicle
  routing problem with route balancing.
\newblock {\em Journal of Heuristics}, 13(5):455--469.

\bibitem[\protect\astroncite{Jozefowiez et~al.}{2008}]{jozefowiez2008}
Jozefowiez, N., Semet, F., and Talbi, E.-G. (2008).
\newblock Multi-objective vehicle routing problems.
\newblock {\em European Journal of Operational Research}, 189(2):293--309.

\bibitem[\protect\astroncite{Jozefowiez et~al.}{2009}]{jozefowiez2009}
Jozefowiez, N., Semet, F., and Talbi, E.-G. (2009).
\newblock An evolutionary algorithm for the vehicle routing problem with route
  balancing.
\newblock {\em European Journal of Operational Research}, 195(3):761--769.

\bibitem[\protect\astroncite{Kovacs et~al.}{2014}]{kovacs2014}
Kovacs, A.~A., Golden, B.~L., Hartl, R.~F., and Parragh, S.~N. (2014).
\newblock Vehicle routing problems in which consistency considerations are
  important: A survey.
\newblock {\em Networks}, 64(3):192--213.

\bibitem[\protect\astroncite{Lacomme et~al.}{2015}]{lacomme2015}
Lacomme, P., Prins, C., Prodhon, C., and Ren, L. (2015).
\newblock A multi-start split-based path relinking ({MSSPR}) approach for the
  vehicle routing problem with route balancing.
\newblock {\em Engineering Applications of Artificial Intelligence},
  38:237--251.

\bibitem[\protect\astroncite{Laumanns et~al.}{2006}]{laumanns2006}
Laumanns, M., Thiele, L., and Zitzler, E. (2006).
\newblock An efficient, adaptive parameter variation scheme for metaheuristics
  based on the epsilon-constraint method.
\newblock {\em European Journal of Operational Research}, 169(3):932--942.

\bibitem[\protect\astroncite{Lin et~al.}{2014}]{lin2014}
Lin, C., Choy, K.~L., Ho, G.~T., Chung, S., and Lam, H. (2014).
\newblock Survey of green vehicle routing problem: past and future trends.
\newblock {\em Expert Systems with Applications}, 41(4):1118--1138.

\bibitem[\protect\astroncite{Matl et~al.}{2016}]{matl2016}
Matl, P., Hartl, R.~F., and Vidal, T. (2016).
\newblock Equity objectives in vehicle routing: A survey and analysis.
\newblock {\em Transportation Science}.
\newblock (forthcoming).

\bibitem[\protect\astroncite{Oyola and L{\o}kketangen}{2014}]{oyola2014}
Oyola, J. and L{\o}kketangen, A. (2014).
\newblock {GRASP-ASP}: An algorithm for the {CVRP} with route balancing.
\newblock {\em Journal of Heuristics}, 20(4):361--382.

\bibitem[\protect\astroncite{Park and Kim}{2010}]{park2010}
Park, J. and Kim, B.-I. (2010).
\newblock The school bus routing problem: A review.
\newblock {\em European Journal of Operational Research}, 202(2):311--319.

\bibitem[\protect\astroncite{Pasia et~al.}{2007a}]{pasia2007a}
Pasia, J.~M., Doerner, K.~F., Hartl, R.~F., and Reimann, M. (2007a).
\newblock A population-based local search for solving a bi-objective vehicle
  routing problem.
\newblock In {\em Evolutionary Computation in Combinatorial Optimization},
  pages 166--175. Springer.

\bibitem[\protect\astroncite{Pasia et~al.}{2007b}]{pasia2007b}
Pasia, J.~M., Doerner, K.~F., Hartl, R.~F., and Reimann, M. (2007b).
\newblock Solving a bi-objective vehicle routing problem by pareto-ant colony
  optimization.
\newblock In {\em Engineering Stochastic Local Search Algorithms. Designing,
  Implementing and Analyzing Effective Heuristics}, pages 187--191. Springer.

\bibitem[\protect\astroncite{Prins}{2004}]{prins2004}
Prins, C. (2004).
\newblock A simple and effective evolutionary algorithm for the vehicle routing
  problem.
\newblock {\em Computers \& Operations Research}, 31(12):1985--2002.

\bibitem[\protect\astroncite{Prins et~al.}{2009}]{prins2009}
Prins, C., Labadi, N., and Reghioui, M. (2009).
\newblock Tour splitting algorithms for vehicle routing problems.
\newblock {\em International Journal of Production Research}, 47(2):507--535.

\bibitem[\protect\astroncite{Samanlioglu}{2013}]{samanlioglu2013}
Samanlioglu, F. (2013).
\newblock A multi-objective mathematical model for the industrial hazardous
  waste location-routing problem.
\newblock {\em European Journal of Operational Research}, 226(2):332--340.

\bibitem[\protect\astroncite{Toth and Vigo}{2003}]{toth2003}
Toth, P. and Vigo, D. (2003).
\newblock The granular tabu search and its application to the vehicle-routing
  problem.
\newblock {\em Informs Journal on Computing}, 15(4):333--346.

\bibitem[\protect\astroncite{Toth and Vigo}{2014}]{toth2014}
Toth, P. and Vigo, D. (2014).
\newblock {\em Vehicle routing: problems, methods, and applications},
  volume~18.
\newblock Siam.

\bibitem[\protect\astroncite{Vidal}{2016}]{vidal2016}
Vidal, T. (2016).
\newblock Technical note: Split algorithm in {O}(n) for the capacitated vehicle
  routing problem.
\newblock {\em Computers \& Operations Research}, 69:40--47.

\bibitem[\protect\astroncite{Vidal et~al.}{2012}]{vidal2012}
Vidal, T., Crainic, T.~G., Gendreau, M., Lahrichi, N., and Rei, W. (2012).
\newblock A hybrid genetic algorithm for multidepot and periodic vehicle
  routing problems.
\newblock {\em Operations Research}, 60(3):611--624.

\bibitem[\protect\astroncite{Zitzler et~al.}{2007}]{zitzler2007}
Zitzler, E., Brockhoff, D., and Thiele, L. (2007).
\newblock The hypervolume indicator revisited: On the design of
  pareto-compliant indicators via weighted integration.
\newblock In {\em International Conference on Evolutionary Multi-Criterion
  Optimization}, pages 862--876. Springer.

\bibitem[\protect\astroncite{Zitzler et~al.}{2003}]{zitzler2003}
Zitzler, E., Thiele, L., Laumanns, M., Fonseca, C.~M., and Da~Fonseca, V.~G.
  (2003).
\newblock Performance assessment of multiobjective optimizers: an analysis and
  review.
\newblock {\em IEEE transactions on evolutionary computation}, 7(2):117--132.

\end{thebibliography}
\end{spacing}

\end{document}